\magnification1200

\magnification1200
\def \sn{{\smallskip \noindent}}
\def \bn {{\bigskip \noindent}}
\def \la {{\longrightarrow}}
\def \O{{\cal O}}
\def \P{{\bf  P}}

\def \Z{{\bf Z}}

\def \dim{{\rm dim}}
\font\medium=cmbx10 scaled \magstep1
\font\large=cmbx10 scaled \magstep2
\vskip .5cm \noindent

{\centerline {{\large Towards a Mori Theory on Compact K\"ahler Threefolds III}}}
\vskip .5cm \vskip .5cm

{\centerline {{\medium Thomas Peternell}}}

\vskip .5cm \vskip .5cm

\noindent {\medium Introduction}
\vskip .5cm 
\noindent In this note we continue the study of the bimeromorphic geometry of compact K\"ahler threefolds. The final 
aim should be - as in the algebraic case - to construct minimal models for threefolds with non-negative Kodaira dimension, describe
the way how to get a minimal model, to prove abundance for the minimal models, i.e. semi-ampleness of the 
canonical bundle, and finally to construction Fano fibrations on appropriate models of threefolds with
$\kappa = - \infty.$ 
\bn Concerning abundance we show that the canonical divisor of a minimal K\"ahler threefold $X$ is good,
i.e. $\kappa (X) = \nu (X),$ where $\nu(X) $ denotes the numerical Kodaira dimension, i.e. the largest 
number $m$ such that $K_X^m \not \equiv 0;$ unless $X$ is simple and non-Kummer. Here $X$ is {\it simple},
if there is no positive-dimensional subvariety through the general point of $X$ and $X$ is {\it Kummer} if
it has a bimeromorphic model which is the quotient of a torus by a finite group. These simple non-Kummer
varieties are expected
not to exist but this can be only a consequence of a completely developped minimal 
model theory in the K\"ahler
case. By a result of Nakayama [Na87] (Kawamata in the algebraic case), it follows from $\kappa (X)
= \nu (X)$ that 
$K_X$ is semi-ample, i.e. some multiple $mK_X$ is generated by global sections. So abundance
holds on K\"ahler threefolds with the possible exception of simple non-Kummer threefolds. 
\bn Furthermore we prove, using essentially Part 1 [CP96] and Part 2 [Pe96] to this paper that a smooth compact K\"ahler threefold
$X$ with $K_X$ not nef carries a contraction unless possibly $X$ is simple and non-Kummer. The 
main steps in the proof are the following
\sn (1) Construction of some curve $C \subset X$ with $K_X \cdot C < 0. $ We distinguish
the case $\kappa (X) \geq 0$ and $\kappa (X) = - \infty$. In the first case we examine carefully
a member in the linear system $\vert mK_X \vert$ to construct $C,$ in the second we use a
result of a recent joint paper of Campana and the author saying that $X$ is uniruled unless
$X$ is simple. In that sedond case it is immediately clear that we can choose $C$ rational.
\sn (2) Next we make $C$ rational (this step works for all compact K\"ahler threefolds).
Here we construct from $C$ a non-splitting family of irreducible curve and examine its
structure. The reason why this family exists is the deformation lemma of Ein-Koll\'ar: 
a curve $C$ in a smooth threefold with $K_X \cdot C < 0$ moves. 
\sn (3) The last step is the construction of a contraction from a rational curve $C$ with
$K_X \cdot C < 0.$ This was in large parts already done in [CP96] and [Pe96]; here we finish
the study in sect. 1 of this paper.
\bn We summarise the results of this paper in the following two theorems.

\bn {\bf Theorem 1} {\it Let $X$ be a minimal K\"ahler threefold (${\bf Q}-$factorial with
at most terminal singularities). Assume that $X$ is not both simple and non-Kummer. Then
$\kappa (X) = \nu (X),$ hence $K_X$ is semi-ample. In particular, if $\kappa (X) = 0$ (and $X$ not simple non-Kummer),
then $K_X \equiv 0$ and $mK_X = \O_X$ for some positive $m$).}

\bn {\bf Theorem 2}  {\it Let $X$ be a  smooth compact K\"ahler threefold  with $K_X$ not nef.
Then $X$ carries a contraction unless (possibly) $X$ is simple with $\kappa (X) = - \infty.$}

\bn By a contraction we mean a surjective map $\varphi : X \la Y$ to a normal compact variety
with connected fibers such that $-K_X$ is $\varphi-$ample and $b_2(X) = b_2(Y) + 1.$
We would like to have that $Y$ is again a K\"ahler space but at the moment we are still
in trouble if $\varphi$ is the blow-up of a smooth curve in the smooth threefold $Y$ - in
that case $Y$ will be only K\"ahler if $\varphi$ is choosen appropriately, namely the 
ray generated by the curves contracted by $\varphi$ must be extremal in the dual cone
to the K\"ahler cone of $X$ (see [Pe96]). 

\bn There are several problems arising. 
\sn (a) First of all we would like to contruct a curve with $K_X \cdot C < 0$ also in 
the "simple" case if $K_X$ is not nef. 
This would mean that we should construct "directly" - i.e. without using any specific 
information on $X$ some curve $C$ with $K_X \cdot C < 0.$ This requires certainly new techniques
and probably a threefold proof would work in any dimension and also with terminal singularities.
\sn (b) We need to prove the existence of contractions also for ${\bf Q}-$Gorenstein threefolds with terminal
singularities in order to perform the Mori program. The Gorenstein case will probably be the
same as the smooth case, but in the presence of non-Gorenstein singularities there obstructions
to moving curves so that some new arguments are needed.
\sn (c) We must overcome the difficulty with the K\"ahler property of $Y$ in case of a
blow-up $\varphi: X \la Y$ along a smooth curve. 
 
\bn As already observed in [Pe96], one major consequence of this programme would be
that simple K\"ahler threefolds are Kummer, in particular there are no simple threefolds
of negative Kodaira dimension, all those being uniruled.

\bn I want to thank C. Okonek for interesting discussions on contact manifolds.  

\bn \bn
{\medium Preliminaries}
\bn {\bf (0.1) } Let $X$ be an irreducible normal compact complex in class ${\cal C},$
e.g. a normal compact K\"ahler space. Let $L$ be a line bundle on $X.$ Then $L$ is nef
if there is a desingularisation $\pi: \hat X \la X$ such that $\pi^*(L)$ is nef. By [Pe96,4.6]
this is independent on the choice of $\pi$ at least in dimension 3. 
\bn {\bf (0.2)} \ ÊA normal compact K\"ahler threefold ($n-$fold) is minimal if $X$ is 
${\bf Q}$-factorial, has only terminal singularities and $K_X$ (i.e. some $mK_X$) is nef.
\bn {\bf (0.3)} We will often use $C_{n,m}$ for K\"ahler threefolds: if $f: X \la Y$ is a surjective
fiber space with $X$ a smooth compact K\"ahler threefold (so $n = 3$) 
then $\kappa (X) \geq \kappa (Y) +
\kappa (F),$ where $F$ denotes the general fiber ($m = \dim Y$). See [Fu78],[Ka81],[Ue87].
\bn {\bf (0.4)}ÊA compact manifold is {\it simple} if there is no positive dimensional subvarieties
through the general point of $X.$ The only known K\"ahler examples arise from tori. To make
this more precise one says that a compact manifold (or normal compact space) is {\it Kummer}
if $X$ is bimeromorphic to a quotient $T/G$ of a torus by a finite group $G.$ The conjecture
is that simple K\"ahler manifolds are Kummer. A standard reference here is [Fu83]. 
\bn {\bf (0.5)} \ ÊSome further notations: $a(X) $ denotes always the algebraic dimension of $X.$
The irregularity of $X$ is $q(X) = \dim H^1(X,\O_X).$ Finally $N_1(X) \subset H_2(X,{\bf R})$ is the vector space 
generated by the classes of
irreducible curve. Inside $N_1(X)$ we have the closed cone
${\overline {NE}}(X)$
generated by the classes of irreducible curves.

\bn \bn 
{\medium 1. Abundance for K\"ahler threefolds}
\vskip .5cm
\noindent In this section we prove the Abundance Conjecture for (non-simple) K\"ahler threefolds.  
First we show

\bn {\bf 1.1 Theorem} {\it Let $X$ be a minimal K\"ahler threefold with $\kappa (X)  = 0.$
Assume that $X$ is not both simple and non-Kummer. Then $K_X \equiv 0.$}

\bn {\bf Proof.} The assertion being known for projective minimal threefolds by Miyaoka [Mi88] and
Kawamata [Ka92], we shall assume that $X$ is non-projective. Then by [CP98] $X$ has a bimeromorphic
model $X'$ with at most quotient singularities such that there is a finite cover $\tilde X \la X',$
\'etale in codimension 1, with $\tilde X$ a torus or a product $E \times S$ of an elliptic curve
$E$ with a K3-surface $S.$ 
In the algebraic case such a conclusion is of course
false: Calabi-Yau threefolds are simply connected. The reason why the non-algebraic case is somehow more
special than the projective is the existence of 2-forms on non-algebraic threefolds (Kodaira's theorem).
Now back in our specific situation, we conclude that $K_{X'}Ê\equiv 0.$ We claim that $X'$ has only 
canonical singularities (a priori we know only that $(X,0)$ has log terminal singularities).
This is seen as follows. Choose $m > 0$ such that $mK_{X'}
\simeq \O_{X'}.$ Let $\pi: \hat X \la X$ be a desingularisation. Write
$$mK_{\hat X}Ê= \pi^*(mK_{X'})Ê+ \sum a_i E_i $$ 
with $E_i$ the exceptional components of $\pi.$ We need to show that $a_i \geq 0$ for all $i.$ 
If some $a_i < 0,$ then $\pi_*(mK_{\hat X})$ is a proper subsheaf of $mK_X = \O_X;$ hence
$H^0(mK_{\hat X} = 0.$ Passing to a high multiple of $m,$ we deduce $\kappa (\hat X) = - \infty,$
a contradiction. 
So $X'$ has only canonical singularities. Take a partial crepant resolution $\pi: \hat X \la X'$ [Re80].
One can even assure that $\hat X$ is ${\bf Q}-$factorial [Re83],[Ka88]. So
$\hat X$ has only terminal singularities and still $K_{\hat X}Ê\equiv 0.$ Now the bimeromorphic 
meromorphic
map $X \rightharpoonup \hat X$ is an isomorphism in codimension 1 [Ha87,Ko89]. Here we use the
fact that $K_X$ is nef! Thus $K_X \equiv 0.$ 

\bn {\bf 1.2 Remark.} 
In case $q(X) > 0,$ it is much easier to conclude; the argument being independent of [CP98]. Here is
the reasoning in that case.
Let $\alpha: X \la A$ be the Albanese map of $X.$ Since $\kappa (X) = 0,$ $\alpha$ has to be surjective
(by $C_{3,1}$ and $C_{3,2}).$ 
Let $f: X \la Y$ be the Stein factorisation of $\alpha$ (we shall se that actually $\alpha$ has connected fibers). 
Notice that $\kappa (F) = 0$ for the general 
fiber $F$ of $f.$ 
\sn (1) Suppose $\dim A = \dim Y = 3.$ Then $X \la A$ is unramified in codimension 1, in fact, otherwise
$K_X$ would contain the ramification divisor $R$ whose image in $A$ is a nef divisor with 
$\kappa > 0,$
so that $\kappa (R) > 0,$ hence $\kappa (X) > 0.$ Of course, this is well-known, at least in the
algebraic case, see [KV80]. So $\alpha = f$ (by the universal property of $\alpha$), hence $\alpha$ is birational. 
Therefore
$$ K_X = \alpha^*(K_A) + \sum_{i \in I} \lambda_i E_i = \sum \lambda_i E_i,$$
where $E_i$ are the exceptional components of $\alpha$ and $\lambda_i > 0.$ Then $K_X$ being nef, it follows
easily that $I = \emptyset.$ So $K_X = \O_X.$
\sn (2) Next suppose that $\dim Y = 1.$ Hence $Y = A$ is an elliptic curve and again $f = \alpha.$ Since $K_F$
is nef and $\kappa (F) = 0,$ we deduce that $K_F$ is torsion. Choose $m > 0$ such that $mK_X$ is Cartier and
that $h^0(mK_X) = 1.$ In particular $mK_F = \O_F.$ Let $s \in H^0(mK_X)$ be a non-zero section. Then $sÊ\vert F = 0$
or $s\vert F$ has no zeroes. Writing
$$ mK_X = \sum a_i D_i$$
with $a_i > 0,$ we conclude that $D_i \cap F = \emptyset,$ i.e. $\dim D_i = 0.$ On the other hand $\sum a_i D_i$
is $\alpha-$nef, and this is only possible if some multiple $kmK_X$ is a multiple of fibers so that
$$ kmK_X = \alpha^*(L)$$
with a line bundle $L$ on $A.$ Since $\kappa (X)  = 0,$ we have $\kappa (L) = 0,$ therefore $L = \O_A.$ 
\sn (3) Finally, let $\dim Y = 2.$ Again we conclude $f = \alpha$ and $Y = A.$ 
\sn (3.1) $a(A) = 0.$ Letting $s$ and $mK_X = \sum a_i E_i$ be as in (2), we again have
$\dim \alpha (D_i) \leq 1$ (now $F$ is an elliptic curve). Since $A$ does not contain compact curves,
we have $\dim \alpha (D_i) = 0$ for all $i,$ which again contradicts the nefness of $K_X.$ 
\sn (3.2) $a(A) = 1.$ Then we use the algebraic reduction $A \la B$ to an elliptic curve and obtain a
map $g: X \la B.$ Now we argue as in (2).

\sn (3.3) $a(A)  = 2.$ By $C_{3,2}^+$ (cp. [Ue87]) we have ${\rm Var}(f) = 0.$ Hence there exists a finite cover
$\tilde A \la A $ and a meromorphic generically finite map $F \times \tilde A \rightharpoonup X$
with an elliptic curve $F.$ Thus $X$ is algebraic, contradiction. So this subcase cannot happen.

\bn {\bf 1.3 Corollary}  {\it Let $X$ be a minimal K\"ahler threefold. Assume that $X$ is not
both simple and non-Kummer. Then $\kappa (X) = \nu (X)$ and therefore $K_X$ is semi-ample.}

\bn As usual, $\nu (X)$ is the numerical Kodaira dimension of $X,$ i.e. the largest number
$k$ such that $K_X^k \not \equiv 0.$ 

\bn {\bf Proof.} Again we may assume $X$ not projective. By [Na87] it suffices to prove
$\kappa (X) = \nu(X).$  
\sn (1) We have $\kappa (X) \geq 0$ by [CP98].
\sn (2) If $\kappa (X) = 0,$ then $\nu(X) = 0$ by (1.1).
\sn (3) If $\kappa (X) \geq 1,$ we can apply [Ka85,7.3] to obtain $\kappa (X) = \nu (X)$; the arguments there 
remain valid in the K\"ahler case (at least for threefolds). Note that $K_X^3  > 0$
forces $X$ to be projective.

\bn {\bf 1.4 Remark} (1) The case that $X$ is simple and not Kummer remains open. The
first difficulty is that we might have $K_X $ nef but $\kappa (X) = - \infty.$ 
Furthermore, if $\kappa (X) = 0,$ then at least $\nu (X) \ne 1.$ In case $X$ is projective, this is proved
in [Mi88]; his arguments work in the K\"ahler case as well. However it might still be
possible that there is a simple non-Kummer minimal threefold with $\kappa (X) = 0$
and $\nu (X) = 2.$ 
\sn Everything would be settled if one could produce a 1-form on a minimal
simple threefold, possibly after finite cover, \'etale in codimension 1 or 2.
In fact, suppose that $q(X) > 0$ and $X$ smooth for simplicity and later reference. 
Then, $X$ being simple, we must have $q(X) = 3,$ the
Albanese $\alpha : X \la A$ has to be onto the threedimensional torus $A.$ Since $A$ has
no divisors, $\alpha$ is unramified in codimension 1, so that actually, arguing as in 
(1.1), $\alpha$ is bimeromorphic and $X$ is Kummer.
\sn (2) We treat one of the open problems in (1) in a special case. Assume that
$X$ is smooth and that $mK_X = \O_X(D)$
with a {\bf smooth} surface (and $X$ simple, $K_X$ nef). By adunction $K_D$ is nef and 
$K_D^2 = 0.$ Therefore $c_2(D) \geq 0.$ On the other hand, $c_2(X) \cdot D = c_2(D)$
and $c_2(X) \cdot D = c_2(X) \cdot mK_X = - 24m \chi(X,\O_X).$ If $q(X) = 0,$ then
$\chi(X,\O_X) > 0,$ hence $c_2(D) < 0,$ contradiction.

\bn {\bf 1.5 Theorem}  {\it Let $X$ be a smooth compact K\"ahler threefold with $\kappa (X)
= 0.$ Suppose that $X$ is not both simple and non-Kummer. Then $X$ has a minimal model.}

\bn {\bf Proof.} The algebraic case being settled by Mori [Mo88] (and previously by Ueno
and Viehweg if $q(X) > 0,$ see [Vi80]), we may assume that $X$ is not projective. 
Then however such a model has already been constructed in the proof of (1.1) (called
$\hat X$ there).

\bn \bn \bn  
{\medium 2. Existence of contractions}
\bn
{\bf (2.1)} Let $X$ be a smooth compact K\"ahler threefold and $C \subset X$ a rational
curve such that $K_X \cdot C < 0.$ Then $C$ determines a non-splitting family of rational
curves [CP96]. In [CP96] and [Pe96] we proved that in turn a non-splitting family $(C_t)_{t \in T}$
of rational curves defines a 
contraction $\varphi: X \la Y$ except in one case. 
This case is the following, called Case (E) below:
\sn $X$ is non-projective, of course, $K_X \cdot C_t = -1,$ the family $(C_t)$ fills up a non-normal surface $S$, furthermore
$\kappa (X) = - \infty$ and $S \cdot C_t \geq 0.$ In this case we expect that $X$ carries
a "generic" conic bundle structure, in particular $X$ is uniruled and it is proved in [Pe96,3.5]
that, given a generic conic bundle structure on $X,$ then $X$ carries a contraction.
\sn The aim of this section is to settle Case (E), so that we obtain:

\bn {\bf 2.2 Theorem}  {\it Let $X$ be a smooth compact K\"ahler threefold, $C\subset X$ a 
rational curve with $K_X \cdot C < 0.$ Then there exists a surjective morphism $\varphi: X \la Y$ 
to a normal compact complex space $Y$ such that
{\item {(1)} $\varphi$ has connected fibers;
\item {(2)} $-K_X$ is $\varphi-$ample;
\item {(3)} $\rho(X) = \rho(Y)  + 1;$
\item {(4)} $Y$ is K\"ahler except possibly if $Y$ is smooth and $\varphi$ is blow-up
along a smooth curve. In this last case $Y$ is K\"ahler if and only if the class of a
fiber is extremal in the closure of the dual of the K\"ahler cone.} }

\bn {\bf Proof.} As explained in (2.1), we have only to deal with the case (E). So
assume that we are in the situation of (E). 
\sn (1) First we claim that it is sufficient to prove that $X$ is uniruled. In fact,
if $X$ is uniruled, then [CP96,2.10] yields $a(X) \ne 2$ and the proof of ibid.(2.12) shows 
that $a(X) \ne 0.$ Thus $a(X) = 1.$ Now [CP96,2.13] shows that $X$ carries a generic
conic bundle structure (we used in (2.13) the assumption of algebraic approximability
only to conclude uniruledness!). Hence a contraction exists as explained in (2.1). 
\sn (2) In case $X$ is not simple, we can now just apply [CP98] to conclude from
$\kappa (X) = - \infty$ that $X$ is uniruled. However our reasoning below for the "simple" 
case also applies here and makes the proof of (2.2) independent of [CP98].
In order to prove uniruledness we consider the normalisation $\nu : \tilde S \la S.$
By [CP96] $\tilde S$ is a smooth minimal ruled surface, so we have a $\P_1-$bundle
structure $g: \tilde S \la \tilde C$ over a smooth curve $\tilde C.$ Let $N \subset S$ be the non-normal
locus, equipped the conductor ideal and $\tilde N \subset \tilde S$ be the analytic
preimage. 
\sn (4) We shall assume here that there is a multisection $C_1 \subset  \tilde N$ of
$g$ and a fiber $F$ of $g$ such that $g(C_1) = g(F).$ We shall use the standard theory
of ruled surfaces as treated in [Ha77,V.2]. Let $C_0$ be a section of $\tilde S$ with
$C_0^2 $ minimal and let $e = - C_0^2.$ Then we can write
$$ C_1 = aC_0 + bF$$
with $a \geq 1$ (here $=$ means numerical equivalence).
Consider the "cycle map" $\nu_*: N_1(\tilde S) \la N_1(X).$ I claim that for every class $C'
\in {\overline {NE}}(\tilde S) $ there is a positive rational number $\mu$ such that
$$ \nu_*(C') = \mu C_t = \mu \nu_*(F).$$ 
By our assumption $$\nu_*(C_1) = \lambda \nu_*(F)$$ 
with a positive rational number $\lambda.$ Hence 
$$\nu_*(aC_0) = \nu_*(C_1) - b\nu_*(F) = (\lambda - b)\nu_*(F).$$  
Let $\omega$ be a K\"ahler form on $X.$ Then 
$$ \int_{C_0}Ê\nu^*(\omega) = ({{\lambda - b}Ê\over {a}}) \int_{C_t}Ê\omega > 0,$$
hence $\lambda > b$ and our claim (*) holds for $C' = C_0.$ 
If now $e \geq 0,$ then ${\overline {NE}}(\tilde S) = {\bf R}_+(C_0) + {\bf R}_+(F),$
and (*) follows. If however $e < 0,$ then 
$$ {\overline {NE}}(\tilde S) = {\bf R}_+C' + {\bf R}_+F$$
with $C' = C_0 + {e \over {2}}F.$ Thus 
$$ \nu_*(C') = \nu_*(C_0) + {e \over {2}}Ê\nu_*(F) = ({{(\lambda - b)} \over {a}} + {e \over {2}}) 
\nu_*(F).$$ 
Integrating the K\"ahler form $\omega$ again, we see that the coefficient is positive, so we are done
with (*). 
Our conclusion from (*) is the following: 
\sn (**) if $L$ is a line bundle on $S$ such that $L \cdot C_t > 0 $ (resp.
$L \cdot C_t = 0)$ then $\nu^*(L)$ is ample, so $L$ is ample (resp. $L \equiv 0$).
\sn Since $S \cdot C_t \geq 0,$ we deduce that either
{\item {(a)}Ê the normal bundle $N_S$ is ample, or
\item {(b)}Ê$N_S \equiv 0.$}
\sn In (a) we have $S^3 > 0,$ hence $X$ is projective [CP96,2.9].
\sn (b) By [Mo82],[Re94], we have $H^1(S,\O_S) = 0$ since $-K_S$ is ample. In particular 
$$ {\rm Pic}(S) \subset H^2(S,{\bf Z}).$$ 
We claim that 
$$N_S  \simeq \O_S. \eqno (+)$$
In order to prove (+) we will verify that $\tilde N$ consists only of $C_0$ and one fiber $F$
so that in particular
$$ H^1(\tilde N,{\bf Z}) = 0. \eqno (++)$$
Assuming (++) for the moment, we consider the Mayer-Vietoris type sequence
$$  H^1(\tilde N,\Z) \la H^2(S,\Z) \la H^2(\tilde S,\Z) \oplus H^2(N,\Z) \la H^2(\tilde N,\Z).$$ 
Then (++) implies
$$ H^2(S,\Z) \subset H^2(\tilde S,\Z) \oplus H^2(N,\Z) \simeq \Z^3$$
(notice that the reduction of $N$ is an irreducible curve because $\tilde N$ has only two
components). Thus $H^2(S,\Z) $ is torsion free.
Moreover $c_1(N_S) = 0 $ and therefore (+) follows from ${\rm Pic}(S) \subset H^2(S,\Z).$ 
From $H^1(\O_S) = 0$ and $H^2(S,\O_S) = H^0(S,\omega_S) = 0$ (note $\omega_S \cdot C_t = -1$),
it follows
$$ \chi(N_S) = \chi(\O_S) = 1,$$
hence $S$ deforms in a 1-dimensional family whose general member $S_t$ has again negative Kodaira
dimension, i.e. $H^0(\omega^{m}_{S_t}) = 0$ for all positive $m$  and therefore (going to a
desingularisation) we see that $S_t$ is uniruled, hence $X$ is uniruled. 
It therefore remains to prove that $\tilde N$ has only two components in order to finish (4).
Write
$$ \nu^*(\omega^*_S) = \alpha C_0 + \beta F$$ 
and
$$ \tilde N = \gamma C_0 + \delta F.$$
Then $\alpha = 1$ (because $\nu^*(\omega^*_S) \cdot F = 1$) and moreover $\gamma > 0$ (because
$\tilde N$ contains $C_1,$ hence $\gamma \geq a.$  
We are going to use the equation
$$ \omega_{\tilde S}Ê= \nu^*(\omega_S) - \tilde N$$
(see [Mo82]). Since 
$$ \omega^*_{\tilde S}Ê= 2C_0 + (e+2-2g),$$
$g$ the genus of $C,$ we obtain
$$ 2 = \gamma + \alpha \ {\rm and} \ e+2-2g = \beta + \delta.$$
Thus $\gamma = 1.$ Notice that $\delta > 0$ since a fiber and a multi-section (now a section)
are identified. Note also that  $\nu^*(\omega_S^*)$ is ample by (**). Now it is a simple
calculation using [Ha77,V.2] to obtain $g = 0, e \geq 0, \beta = e+1$ and $\delta = 1.$ 
So $\tilde N = C_0 + F,$ even ideal-theoretically. This proves (++).
\sn (5) Finally we have to treat the case that no multisection of $g$ is identified with a fiber
by $\nu.$ Then $g$ induces a map $h: S \la B$ to a - usually non-normal - curve $B$ 
and a map $\mu: C \la B,$ the normalisation of $C,$ such that
$$h \circ \nu = \mu \circ g.$$
So for every $t \in T$ we find some $b \in B$ such that $C_t \subset S_b:= h^{-1}(b),$
and $S_b$ consists only of $C_t's.$ We fix a general smooth point $b \in B.$
Then $S_b,$ the analytic fiber over $b,$ is of the form 
$$ S_b = C_{t_1} + \ldots C_{t_r},$$
and of course $S_b$ is as a Cartier divisor in a Gorenstein surface free from embedded points.
Actually $S_b$ is locally a complete intersection (in $X).$ 
\sn (a) If $S_b$ is reducible, then
$$ -K_X \cdot S_b \geq 2,$$
hence by Ein's deformation lemma (see [Ko96,II.1.16]) the cycle $S_b$ of rational curves
deforms in an at least 
2-dimensional family. Since the deformations of $S_b$ inside $S$ have dimension 1,
we conclude the existence of a covering family $(B_s)$ of curves for $X$ (the general
$B_s$ might however not be rational). So $X$ is
not simple, the only case to be excluded if we use [CP98]. Here is the way to avoid
[CP98]: let $u: \tilde S_b = \nu^{-1}(S_b) \la S_b$ be the normalisation of $S_b.$
Then $\tilde S_b$ consists of $r$ disjoint smooth rational curves. Now 
$$ \dim_{[u]} {\rm Hom}(\tilde S_b,X) \geq -K_X \cdot S_b + 3 \chi(\O_{\tilde S_b})
= -K_X \cdot \tilde S_b + 3r,$$
see e.g. [Ko96]. Taking into account $\dim {\rm Aut}(\tilde S_b) = 3r,$ 
we can still conclude the existence of a covering family of reducible rational curves
for $X$ so that $X$ is uniruled. 
\sn (b) Suppose finally that $S_b$ is irreducible. Having in mind that $S_b$ is 
Cartier in $S$ and also locally a complete intersection, we consider the exact
sequence of conormal bundles
$$ 0 \la N^*_{S \vert X}Ê\vert S_b \la N^*_{S_b \vert X}Ê\la N^*_{S_b \vert S}Ê\la 0.$$
Since $N^*_{S_b \vert S} = \O_{S_b}$ and since $N^*_{S \vert X}Ê\vert S_b$ has nef dual,
it follows via the sequence that $N_{S_b \vert X}$ is nef. This contradicts easily
$-K_X \cdot C_t = -1.$ 

\bn \bn \bn 
{\medium 3. Existence of rational curves}

\bn 
In this section we prove that any smooth compact K\"ahler threefold $X$ whose canonical bundle is
not nef, carries a rational curve $C$ with $K_X \cdot C < 0,$ unless $X$ is simple with negative
Kodaira dimension. We first treat the case $\kappa (X) \geq 0.$ 

\bn {\bf 3.1 Theorem} {\it Let $X$ be a smooth compact K\"ahler threefold with $\kappa (X) \geq 0.$
If $K_X $ is not nef, there exists a rational curve $C \subset X$ with $K_X \cdot C < 0,$ and therefore
there is a contraction $\varphi: X \la Y$.}

\bn {\bf Proof.} Since $\kappa (X) \geq 0,$ we find $m \in {\bf N}$ Êsuch that
$$ mK_X = \O_X(\sum \lambda_i D_i),$$
where $\lambda_i \in {\bf N}$ and $D_i \subset X$ are irreducible. By [Pe96,4.9] there exists some 
irreducible component   $D_{i_0}$ such that $K_X \vert D_{i_0}$ is not nef. By renumbering we may assume
that the components with this property are exactly $D_1, \ldots, D_r.$ 
Fix some $1 \leq i \leq r$ and put $D = D_i.$ Let $\nu: \tilde D \la D$ be the 
normalisation and $\pi: \hat D \la \tilde D $ the minimal desingularisation. From adjunction we immediately
obtain that
$$ K_D = r K_X \vert D - E, \eqno (*)$$
where $r > 1 $ is a positive rational number and $E$ an effective ${\bf Q}-$divisor. 
Next observe that
$$ K_{\hat D}Ê= \pi^* \nu^*(K_D) - E'$$
with another effective divisor, so that, putting $L = \pi^* \nu^*(K_X \vert D),$
we obtain a formula of ${\bf Q}-$divisors
$$ L = aK_{\hat D}Ê+ A \eqno (**)$$
with a positive rational number $a < 1$ and an effective ${\bf Q}-$divisor $A.$ 
\sn (1) In a first step we construct an irreducible curve $C \subset D$ such that $K_X \cdot C < 0.$
If $a(D) = 2,$ the existence of $C$ is clear, $K_X \vert D$ being not nef. 
\sn (1.1) Suppose next $a(D) = 0.$
Then $\hat D$ is bimeromorphically a K3-surface or a torus. Then by (**), $L$ has to be effective :
there are curves $C_i \subset \hat D$ and  positive integers $p$ and $m_i$ such that
$$pL = \sum m_i C_i.$$
Observe that every curve in $\hat D$ is rational! So by [Pe96,4.9] we obtain a rational curve
$C \subset \hat D$ with $L \cdot C < 0, $ hence $K_X \cdot \nu(\pi(C)) < 0. $ Thus in case $a(D) = 0$
for some component $D$ we are already settled with the theorem. 
\sn (1.2) Now suppose $a(D) = 1$ and let $f: \hat D \la B$ be the algebraic reduction, an elliptic
fibration. Since $L$ is not nef, [Pe96,4.13] implies that either there exists an irreducible curve
$C \subset \hat D$ with $L \cdot C < 0$ or $L = f^*(G)$ with $G^*$ ample on $B.$ This second alternative
is clearly impossible since $L$ is effective. Hence (1) holds also in this case.

\sn (2)  In all what follows we may assume $C$ irrational, otherwise we are done.
Now $C$ deforms in an at least 1-dimensional family in $X,$ therefore we can "extract"
a  maximal non-splitting family $(C_t)$ of irreducible curves such that $K_X \cdot C_t < 0.$ Here 
"maximal" means that no deformation of the general $C_t$ splits. In other words, we can choose $T$ to be (or
rather project to) an
irreducible component of the cycle space (which is automatically compact, $X$ being K\"ahler). 
By $K_X \cdot C_t < 0$ all the $C_t$ have to be in one of the $D_j's$ for $1 \leq j \leq r,$ 
say in $D.$ By (**), 
we have $K_D \cdot C_t < 0.$
Let $\hat C_t$ be the strict transform of a general $C_t$ in $\hat D.$ Then 
$$ L \cdot \hat C_t < 0.$$ 
In particular (**) yields $K_{\hat D} \cdot \hat C_t < 0.$ $\hat C_t$ being irrational, $\hat C_t$
deforms in an at least 1-dimensional family, see [Ko96,II.1.14,1.15], and we can either 
extract a maximal non-splitting family
$(B_s)$ of irrational curves on $\hat D$ such that $$L \cdot B_s < 0$$ or we obtain the rational curve in $D$ we are looking for. 
To be more precise, assume that we have a splitting.
Take a splitting component $C'$ with $L \cdot C' < 0.$ If $C'$ is irrational, then $C'$ deforms in
a family and we  proceed as before (inductively). If $C'$ is rational, then its image $C''$ in $D$
is a rational curve with $K_X \cdot C'' < 0$ are we are done. Notice however that on surfaces $X$
it is sometimes not possible to extract a non-splitting family from a splitting family $(C_t)$ with
$K_X \cdot C_t < 0.$ The reason is that there are curves with $K_X \cdot C < 0$ which do not move,
namely $(-1)-$curves, and these are the only ones. 
 
Now we pass to a minimal model $\sigma: \hat D \la D_0$ and consider the induced family $(B'_s)$. 
Clearly $(B'_s)$ is a non-splitting family. 
Suppose first that $\hat D$ itself is not minimal, i.e. $\sigma$ is not an isomorphism. Let $p \in D_0$
be a point blown up by $\sigma.$ Then $p \in B'_s$ for some $s$. Since $(B_s)$ is non-splitting, 
we must have $p \in B'_s$ for all $s.$ Then Lemma 3.3 gives a contradiction. 

So $\hat D$ must be minimal. If $\hat D $ is the projective plane, then clearly $(B_s)$ has
to be a family of lines. Anyway, then $\hat D = \tilde D$ and $L^*$ is ample in that case, so
the existence of $C$ rational with $K_X \cdot C < 0$ is clear. So $\hat D$ is ruled with ruling
$p: \hat D \la B.$ By Lemma 3.4, its invariant $e < 0,$ in particular $\hat D$ has no exceptional 
curves and $\tilde D = \hat D.$ Let $F$ be a fiber of $p$ and $C_0$ be a section with $C_0^2 = - e$
minimal. Let $Z$ denote a class in the boundary component of
${\overline {NE}}(\hat D)$ different from the ray generated by $F.$ Then after rescaling
$$ Z \equiv C_0 + {e \over {2}}F.$$
Since $K_{\hat D}Ê\cdot Z \geq 0,$ we conclude from (**) easily that
$$ L \cdot Z \geq 0$$
(consider curves near to the ray ${\bf R}_+(Z).)$ Using our assumption $L \cdot F \geq 0$
(otherwise we are done), we conclude that $L $ is nef, contradicting $L \cdot \hat C_t < 0.$

\bn Notice that in the second part of the proof of 3.1, where we deduced the existence
of rational curve $C\subset X$ with $K_X \cdot C < 0$ from the existence of some curve
$C $ with $K_X \cdot C < 0$ we did not use the assumption $\kappa (X) \geq 0.$
Therefore we have
\bn {\bf 3.2 Theorem}  {\it Let $X$ be a smooth compact K\"ahler threefold. Assume that there
is an irreducible curve $C \subset X$ with $K_X \cdot C < 0.$ Then there is a rational
curve $C \subset X$ with $K_X \cdot C < 0$ and hence there is a contraction on $X.$}  

\bn We still have to prove the following two lemmata.

\bn {\bf 3.3 Lemma} {\it Let $S$ be a ruled surface, $p \in S.$ Let $(C_t)_{t \in T}$ be family of
curves with $C_t$ irreducible for general $t$ and $T$ compact (as always). Assume that $p \in C_t$
for all $t$. Then $(C_t)$ splits. }

\bn {\bf Proof.} Assume that $(C_t)$ does not split. Of course we may restrict ourselves to the
case $\dim T = 1.$ Let $f: S \la B$ be a ruling, $F_b$ the fiber over $b \in B.$ Since
$C_t$ cannot be a fiber, it must be a multi-section of $f$ of degree $d$ and therefore $C_t 
\cap F_b$ consists of $d$ points, counted with multiplicity. Denoting by ${\rm Hilb}_d(X)$
the Hilbert scheme of $d$ points in a projective manifold $X,$ we obtain for every $b$ a map 
$$ \varphi_b: T \la {\rm Hilb}_d(F_b).$$
Let $b_0 = f(p).$ Choose an open neighborhood $U$ of $b_0$ in $B$ and a trivialisation
$S_U \simeq U \times \P_1.$ Then we obtain a map
$$ \Phi: U \times T \la {\rm Hilb}_d(\P_1)$$
such that (by abuse of notation) $\Phi \vert \{b\}Ê\times T = \varphi_b.$ 
Now notice that ${\rm Hilb}_d(\P_1) = S^d\P_1,$ the $d-$th symmetric product of $\P_1.$
Since all $C_t$ pass through $p,$ we conclude
$$ \varphi_{b_0}(T) \subset \{p\}Ê\times S^{d-1}\P_1 \subset S^d\P_1.$$
Using the identification $S^d\P_1 = \P_d,$ the curve $\varphi_{b_0}(T)$ is contained
in a hyperplane, i.e. degenerates. After shrinking $U$, there is no point $q \in F_b, b \in U, b \ne b_0$
such that all $C_t$ pass through $q.$ Therefore $\varphi_b(T)$ is not degenerate for
$b \in U, b \ne b_0.$ This leads to a contradiction since the $\varphi_b(T)$ form a
family of curves: first if $d = 2,$ then $\varphi_{b_0}(T)$ is a line, whereas 
$\varphi_b(T)$ is not for $b \ne b_0;$ if $d \geq 3,$ choose a generic projection
$g: \P_d \rightharpoonup \P_2$ and consider the induced holomorphic maps $g \circ
\varphi_b.$ 

\bn {\bf 3.4 Lemma}  {\it Let $S$ be a ruled surface, $(C_t)_{t \in T}$ be a maximal (i.e. $T$
projects to an irreducible component of the cycle space) non-splitting family
of curves not consisting of fibers such that $K_S \cdot C_t < 0.$ Then $S$ comes from a 
stable vector bundle, i.e. its invariant $e$ is negative.}

\bn {\bf Proof.}    We write for numerical equivalence
$$ C_t \equiv kC_0 + bF$$
where $C_0$ is a section with minimal $C_0^2$ and $F$ is a fiber of $f,$ a ruling. 
Again let $e = - C_0^2.$
More precisely we write
$$ \O_S(C_t) = \O_S(C_0) \otimes f^*(G_t)$$
with a line bundle $G_t$ of degree $b$ on $B.$
By virtue of
$$ K_S = -2 C_0 + (2g-2-e)F $$
we obtain
$$ K_S \cdot C_t = 2ke + k(2g-2-e) - 2b.$$ 
Hence $K_S \cdot C_t < 0$ translates into
$$ b > {{ke}Ê\over {2}}Ê+ k(g-1). \eqno (*) $$ 
If now $e \geq 0,$ i.e. $E$ is not stable, then (*) implies $b > g-1,$ hence
Riemann-Roch gives $h^0(G_t) > 0$ and therefore we obtain a splitting inside the family.

\bn We now turn to the case $\kappa (X) = - \infty.$ 

\bn {\bf 3.5 Theorem}  {\it Let $X$ be a smooth compact K\"ahler threefold or a normal compact 
K\"ahler threefold which is ${\bf Q}$-factorial with at most terminal singularities.
Assume $\kappa (X) = - \infty$ and $X$ not simple. Then there exists a rational curve
$C \subset X$ with $K_X \cdot C < 0.$}

\bn {\bf Proof.}  Of course we have only to treat the non-projective case. By [CP98], $X$
is uniruled. Hence our claim is clear, taking $C$ to be a general member of a covering
family of ratinal curves of $X.$

\bn {\bf 3.6 Remarks} (1) In order to complete the picture that a compact K\"ahler threefold
$X$ with $K_X$ not nef carries a rational curve $C$ with $K_X \cdot C < 0$ it remains to
prove that a simple K\"ahler threefold with $\kappa (X) = - \infty$ carries some curve
(not necessarily rational) with $K_X \cdot C < 0.$ Of course this is a strange situation,
which is expected not to exist. In (2) and (3) we treat some special cases.
\sn (2) There is no simple compact K\"ahler manifold $X$ with $\kappa (X) = - \infty$ such that
$-K_X$ is hermitian semi-positive. In fact if $X$ is simple, then $q(X) = 0 $ and the 
same holds for every finite \'etale cover of $X$, hence [DPS96] implies $H^2(X,\O_X) = 0,$
contradiction. 
\sn (3) The same should be true if $-K_X$ is only nef. We may clearly assume $q(X) = 0.$ Since 
$$ \chi(X,-mK_X) = (2m+1)\chi(X,\O_X) $$
we obtain a contradiction if we can prove that
$$ h^2(-mK_X) \leq C \eqno (*)$$
for a positive constant $C$ and all positive integers $m$, at least if $m$ is 
sufficiently divisible.

\bn We finally summarize the results of this section:

\bn {\bf 3.7 Theorem}  {\it Let $X$ be a smooth compact K\"ahler threefold with $K_X$
not nef. Then there exists a rational curve $C \subset X$ with $K_X \cdot C < 0$ and
therefore $X$ carries a contraction unless (possibly) $X$ is simple with $\kappa (X) =
-\infty.$}     

\bn \bn {\medium 4. K\"ahler contact threefolds}

\bn Let $X$ be a smooth compact K\"ahler threefold. We recall the notion of a contact
structure on $X$ and refer to [Le95] for details. Assume that there is a line bundle $L$ with
$-K_X = 2L.$ A contact structure on $X$ is an $L-$valued 1-form 
$$ \omega \in H^0(\Omega^1_X \otimes L)$$
without zeroes such that locally $\omega \wedge d\omega $ is nowhere $0.$ We therefore have a vector bundle $F$ of rank 2 and an exact sequence
$$ 0 \la F \la T_X \la L \la 0.$$  
Ye [Ye94] proved that a projective contact threefold is either $\P_3$ or the projectivised
tangent bundle of any projective surface. We extend this results to the K\"ahler case;
however we cannot handle at the moment the possibility of a simple contact threefold with
negative Kodaira dimension.

\bn {\bf 4.1 Theorem}  {\it Let $X$ be a smooth compact K\"ahler threefold admitting a
contact structure. Then $X$ is one of the following
{\item {(1)} $X = \P_3$;
\item {(2)} $X$ is the projectivised tangent bundle of a compact K\"ahler surface;
\item {(3)}Ê$X$ is simple with $\kappa (X) = - \infty.$}}

\bn Of course case (3) should not exist. Note also that the projectivised tangent bundle
of a K\"ahler surface has always a contact structure, see e.g. [Le94].
 
\bn {\bf Proof.} The projective case being proved in [Ye95], we shall assume $X$ not 
projective. We may also assume that $X$ is not simple with $\kappa (X) = - \infty.$
\sn (1) First notice that the arguments in [Ye95] show that $\kappa (X) \leq 0.$
In fact, the necessary results of Bogomolov are true in the K\"ahler case, too.
\sn (2) Suppose that $K_X$ is not nef. Then by (3.7) we find a contraction 
$\varphi: X \la Y.$ Since $K_X$ is divisible by 2, $\varphi$ must be the blow-up of
a smooth point or (analytically) a $\P_1-$bundle over a K\"ahler surface. 
Assume first that $\varphi$ is the blow-up of a point with exceptional divisor $E.$
Using the notations from the beginning of the section, we have an exact sequence
$$ 0 \la N^*_E \otimes L \la \Omega^1_X \otimes L \vert E \la 
\Omega^1_E \otimes L \vert E \la  0.$$
Now $L \vert E = \O_E(1),$ hence 
$$ H^0(\Omega^1_X \otimes L \vert E) = H^0(N^*_E \otimes L).$$
But since $N^*_E \otimes L = \O (2),$ every section in $\Omega^1_X \otimes L$ must
have zeroes, contradiction. Hence $\varphi$ is a $\P_1-$bundle over a K\"ahler surface
$Y$ and it follows as in [Ye95] that $X$ is actually $\P(T_Y).$ 
\sn (3) It remains to treat the case that $K_X$ is nef and $\kappa (X) = 0.$ If $X$ is not
simple, then $K_X \equiv 0$ by (1.1). Now a finite \'etale cover of $X$ is either a torus
or a product of a K3-surface with an elliptic curve which immediately gives a contradiction.
So we will suppose that $X$ is simple. 
The contact form is a section in $H^0(X,\Omega^1_X \otimes L)$ with $-K_X = 2L.$
Hence we have an injection $L^* \la \Omega^1_X$. If we knew that $H^0(L^*) \ne 0$
then $q(X) \geq 1$ and $X$ is a torus (see (1.1),(1.4)). However in general we know
only that $\kappa (L^*) = \kappa (X) = 0,$ i.e. $H^0(L^{*m}) \ne 0$ for some positive
$m.$ Then we take a generically finite map $f: \tilde X \la X$ such that 
$$H^0(\tilde X, f^*(L^*)) \ne 0$$ and $\tilde X$ is again smooth.
Since $f^*(L^*) \subset f^*(\Omega^1_X)$ we obtain via the canonical injective map
$f^*(\Omega^1_X) \la \Omega^1_{\tilde X}$ an injection
$$ f^*(L^*) \subset \Omega^1_{\tilde X}$$
and therefore $q(\tilde X) > 0.$ Now $\tilde X$ must again be simple, in particular
$\kappa (\tilde X) = 0.$ Hence $\tilde X$ is Kummer, more precisely the Albanese map
$\alpha: \tilde X \la A$ is bimeromorphic onto the threedimensional torus $A.$ Therefore
$K_{\tilde X}Ê= B$ with $B$ an effective divisor supported on the exceptional locus of
$\alpha$. On the other hand $$K_{\tilde X}Ê= f^*(K_X) + R$$
with another effective divisor $R.$ Since $K_X$ is nef, it follows easily that $R = B$
and thus $f^*(K_X) = \O_X,$ i.e. $K_X \equiv 0,$ or directly a covering $A \la X.$ 
So $X$ is Kummer.

\bn {\bf 4.2 Remark} Assume that $X$ is a simple threefold, $\kappa (X) = - \infty$ and
$X$ admits a contact structure. If there is a curve $C$ with $K_X \cdot C < 0,$
then by (3.2) we have a contraction $\varphi: X \la Y$. This is ruled out as in the
proof of (4.1). Hence we can at least say that $K_X \cdot C \geq 0$ for all curves
$C \subset X.$  

\bn {\bf 4.3 Remark}  Observe that we have proved in (4.1) the following slightly more general 
statement.
\sn Let $X$ be a smooth compact K\"ahler threefold, $L \in {\rm Pic}(X)$ such that $-K_X = 2L.$
Let $\omega \in H^0(\Omega^1_X \otimes L)$ be a section without zeroes. Then $X$ is one
the following
{\item {(a)}Ê$X = \P_3;$
\item {(b)}Ê$X = \P(T_S),$ where $S$ is a K\"ahler surface;
\item {(c)}Ê$X$ is an \'etale "undercover" of a torus or a product of an elliptic curve
with a K3-surface;
\item {(d)} $X$ is simple}

\vfill \eject {\medium References}

\bn
 
\item {[CP96]} Campana,F;Peternell,T.: Towards a Mori theory on compact K\"ahler threefolds,I.
 Math. Nachr. 187, 29-59 (1997)

\item {[CP98]} {\obeylines Campana,F.;Peternell,T.: Holomorphic 2-forms on complex threefolds.
To appear in J. Alg. Geom. 1999}

\item {[DPS94]} Demailly,J.P.;Peternell,T.;Schneider,M.: Compact complex manifolds with
numerically effective tangent bundles. J. Alg. Geom. 3, 295 - 345 (1994)

\item {[DPS96]} Demailly,J.P.;Peternell,T.;Schneider,M.: Compact K\"ahler manifolds with hermitian
semipositive anticanonical bundle. Comp.math. 101, 217-224 (1996) 

\item {[Fu83]} Fujiki,A.: On the structure of compact complex manifolds in class ${\cal C}.$
Adv. Stud. Pure Math. 1, 231 - 302 (1983)

\item {[Fj78]}ÊFujita,T.: K\"ahler fiber spaces over curves. J. Math. Soc. Japan 30, 779 - 794
(1978)

\item {[Ha77]} Hartshorne,R. Algebraic geometry. Graduate Texts in Math. 52. Springer 1977

\item {[Ha87]} Hanamura,M.: On the birational automorphism groups of algebraic varieties.
Comp. math. 63, 123-142 (1987)

\item {[Ka81]} Kawamata,Y.: Characterisation of abelian varieties. Comp. math. 43, 253 - 276
(1981)

\item {[Ka85]}ÊKawamata,Y.: Pluricanonical systems on minimal algebraic varieties. Inv. math.
79, 567-588 (1985)

\item {[Ka85a]} Kawamata,Y.: Minimal models and the Kodaira dimension of algebraic
fiber spaces. J. reine u. angew. Math. 363, 1 - 46 (1985)

\item {[Ka88]} Kawamata,Y.: Crepant blowing ups of threedimensional canonical singularities
and applications to degenerations of surfaces. Ann. Math. 119, 603-633 (1988)
 
\item {[Ka92]} Kawamata,Y.: Abundance theorem for minimal threefolds. Inv. math. 108,
229 - 246 (1992)

\item {[KMM87]} Kawamata,Y.;Matsuda,K.;Matsuki,K.: Introduction to the minimal model problem.
Adv. Stud. Pure Math. 10, 283 - 360 (1987)

\item {[KV80]} Kawamata,Y.;Viehweg,E.: On a characterisation of an abelian variety in the
classification theory of algebraic varieties. Comp.math. 41, 355-359 (1980)

\item {[Ko89]} Koll\'ar,J.: Flops. Nagoya Math. J. 113, 15-36 (1989)

\item {[Ko96]}ÊKoll\'ar,J.: Rational curves on algebraic varieties. Erg. d. Math. vol 32.
Springer 1996 

\item {[Le95]} LeBrun,C.: Fano manifolds, contact structures and quaternionic geometry.
Int. J. Math. 6,419-437 (1995)

\item {[Mi87]} Miyaoka,Y.: The Chern classes and Kodaira dimension of a minimal
variety. Adv. Stud. Pire Math. 10, 449 - 476 (1987)  
Adv. Stud. Pure Math. 10, 

\item {[Mi88]} Miyaoka,Y.: Abundance conjecture for threefolds: $\nu = 1$ case. Comp. math. 68,
203 - 220 (1988)

\item {[Mo82]} Mori,S.: Threefolds whose canonical bundles are not numerically effective.
Ann. Math. 116, 133 - 176 (1982)

\item {[Mo88]} Mori,S.: Flip theorem and the existence of minimal models for 3-folds.
J. Amer. Math. Soc. 1, 117 - 253 (1988)

\item {[Na87]} Nakayama,N.: The lower semi-continuity of the plurigenera of complex
varieties. Adv. Stud. Math. 10, 551-590 (1987) 

\item {[Pe96]}ÊPeternell,T.: Towards a Mori theory on compact K\"ahler threefolds II.
Math. Ann. 311, 729-764 (1998)

\item {[Re80]} Reid,M.: Canonical threefolds. In: Geometrie alg\'ebrique Angers 1980.
Sijthoff and Noordhoff, 273-310 

\item {[Re83]} Reid,M.: Minimal modles of threefolds. Adv. Stud. Pure Math. 1, 131-180 (1983)
                           
\item {[Re94]} Reid,M.: Singular del Pezzo surfaces. Publ. RIMS   , 695-728 (1994) 

\item {[Ue87]} Ueno,K.: On compact analytic threefolds with non-trivial Albanese torus.
Math. Ann. 278, 41 - 70 (1987)

\item {[Vi80]} Viehweg,E.: Klassifikationstheorie algebraischer Variet\"aten der Dimension 3.
Comp. math. 41, 361-400 (1980)

\bn \bn \bn 

\noindent Thomas Peternell \sn Mathematisches Institut  \sn der Universit\"at Bayreuth, 
\sn D-95440 Bayreuth, Germany 
\sn email: thomas.peternell@uni-bayreuth.de

\end